\documentclass[12pt]{amsart}

\usepackage{amsthm}
\usepackage{amssymb}
\usepackage{amsmath}
\numberwithin{equation}{section}
\newcommand{\bea}{\begin{eqnarray}}
\newcommand{\eea}{\end{eqnarray}}
\newcommand{\be}{\begin{eqnarray*}}
\newcommand{\ee}{\end{eqnarray*}}
\newtheorem{theorem}{Theorem}[section]
\newtheorem{lemma}{Lemma}[section]

\newtheorem{proposition}{Proposition}[section]

\begin{document}
\title[Coincidence Indices]{Indices of Coincidence Isometries of the Hyper Cubic Lattice $\mathbb{Z}^{n}$}
\author[Yi Ming Zou]{Yi Ming Zou}
\address{Department of Mathematical Sciences, University of Wisconsin, Milwaukee, WI 53201, USA} \email{ymzou@uwm.edu}
\maketitle

\begin{abstract}
The problem of computing the index of a coincidence isometry of the hyper cubic lattice $\mathbb{Z}^{n}$ is considered. The normal form of a rational orthogonal matrix is analyzed in detail, and explicit formulas for the index of certain coincidence isometries of $\mathbb{Z}^{n}$ are obtained. These formulas generalize the known results for $n\le 4$.
\end{abstract}

\section{Introduction}
\par
The theory of coincidence site lattice (CSL) can be used to describe certain phenomena that arise in the physics of interfaces and grain boundaries (see Bollmann (1970) and Grimmer (1973)). Mathematically, CSL theory concerns the relationship between a lattice $L$ and a transformed copy $\mathcal{A}L$ of $L$, where $\mathcal{A}$ is a linear transformation of the $n$-dimensional real vector space $V$ spanned by $L$. We call $\mathcal{A}$ a {\it coincidence symmetry} if $\mathcal{A}$ is an automorphism of $V$ and $L\cap \mathcal{A}L$ is a sublattice of $L$ with finite index. It is known (see section 2 below) that $\mathcal{A}$ is a coincidence symmetry if and only if the matrix $A$ of $\mathcal{A}$ under a basis of $L$ is a rational matrix. The set of all coincidence symmetries (or the set of all $n\times n$ coincidence matrices) of $L$ forms a group under the multiplication defined by composition (or the multiplication of matrices). If $L$ is a lattice of the Euclidean space $\mathbb{R}^{n}$, then one is interested in the isometries of $\mathbb{R}^{n}$ which are coincidence symmetries of $L$. In this case, we have the {\it coincidence isometry} subgroup formed by all the coincidence isometries (Baake (1997)). 
\par
One of the main problems in CSL theory is the computation of the index of coincidence of $L\cap\mathcal{A}L$ in $L$ (also called degree). M. A. Fortes (1983) provided a general approach to this problem by using the normal form of an integer matrix, and Duneau {\it et al.} (1992) gave a further study along a similar line. Although theoretically it is possible to compute the index of a coincidence transformation via the normal form of the corresponding integer matrix by Fortes' result, no general index formula in $n$-dimension is known even for the coincidence symmetries of the hyper cubic lattice $\mathbb{Z}^{n}$. For the coincidence isometries, it is possible to give more explicit results. Pleasants {\it et al.} (1996) used number theory to treat the planar case. Baake (1997) provided a solution to this problem for the coincidence isometries for dimensions up to 4 by using the factorization properties of certain number systems. However, the method does not generalize to higher dimensions. Recently, geometric algebra method was introduced into the study of CSL theory by Arag\'{o}n {\it et al.} (2001) and Rodriguez {\it et al.} (2005). But only the planer case was treated. Zeiner (2006) provided a detailed analysis for the coincidence indices of hypercubic lattices in 4 dimensions.  
\par
In this paper, we derive several formulas for the index of a coincidence isometry of the lattice $\mathbb{Z}^{n}$ for arbitrary $n$. We analyze the normal form of the corresponding integer matrix of a coincidence isometry taking into account of the orthogonal property. The main formulas are given in Theorem 3.1 and 3.2. 
\par
\section{Preliminaries}
The set of real numbers (respectively, integers and rational numbers) is denoted by $\mathbb{R}$ (respectively, $\mathbb{Z}$ and $\mathbb{Q}$), the set of all non-singular $n\times n$ real matrices is denoted by $GL_{n}(\mathbb{R})$, and the set of $n\times n$ real orthogonal matrices is denoted by $ O_{n}(\mathbb{R})$. Notation for matrices over $\mathbb{Q}$ and $\mathbb{Z}$ are defined similarly. For a nonzero integer matrix $Z$, we denote by $\gcd(Z)$ the greatest common divisor of the nonzero entries of $Z$.
\par
By an $n$-dimensional lattice $L$ with basis $(a_{1}, \ldots, a_{n})$, we mean the free abelian group $\oplus_{i=1}^{n}\mathbb{Z}a_{i}$. In this paper, we only consider lattices in the $n$-dimensional Euclidean space $\mathbb{R}^{n}$, and we assume the lattices are also $n$-dimensional. Thus, a lattice $L\subset \mathbb{R}^{n}$ is given by an $n\times n$ non-singular matrix $A$ (called the structure matrix of $L$), and a basis of the lattice is  
\bea
(a_{1}, \ldots, a_{n})= (e_{1}, \ldots, e_{n})A, 
\eea 
where $(e_{1},\ldots,e_{n})$ is the canonical basis of $\mathbb{R}^{n}$.
\par
By a sublattice $L^{\prime}\subset L$ we mean a subgroup $L^{\prime}$ of finite index in the abelian group $L$. The CSL theory concerns the problems that arise when the intersection $L_{1}\cap L_{2}$ of two lattices happens to be a sublattice of both lattices $L_{1}$ and $L_{2}$. If this is the case, we say that $L_{1}$ and $L_{2}$ are {\it commensurate lattices}. 
\par
Suppose that $L_{i}$ is given by the structure matrix $A_{i}$ ($i=1,2$), and let the basis of $L_{i}$ be $\mathbf{B}_{i}$, i.e. 
\be
\mathbf{B}_{i}=(e_{1}, \ldots, e_{n})A_{i}, \quad i=1,2. 
\ee
Then a theorem due to Grimmer states that $L_{1}$ and $L_{2}$ are commensurate if and only if $A_{2}^{-1}A_{1}$ is a rational matrix. This implies that if $L$ is a lattice with basis $(a_{1},\ldots, a_{n})$ and $A$ is an $n\times n$ non-singular real matrix, then the lattice with basis $(a_{1},\ldots, a_{n})A$ and the lattice $L$ are commensurate if and only if $A$ is a rational matrix.
\par
Consider a lattice $L$ in $\mathbb{R}^{n}$ with the structure matrix $A$. Let $\mathcal T$ be a linear transformation of $\mathbb{R}^{n}$, and let $T$ be the matrix of $\mathcal T$ under the canonical basis $(e_{1},\ldots, e_{n})$. Then the structure matrix of the lattice $\mathcal{T}(L)$ is $TA$. Thus the lattice $\mathcal{T}(L)$ and the lattice $L$ are commensurate if and only if $A^{-1}TA$ is rational. The isometries of $\mathbb{R}^{n}$ that provide commensurate lattices to a lattice $L$ are of special interest (cf. Baake (1997), Arag\'{o}n {\it et al.} (2001), and Rodriguez {\it et al.} (2005)), let us recall the relevant definitions. The following group was defined in Baake (1997):
\be
OC(L)=\{Y\in O(n): [L: L\cap YL]<\infty\}.
\ee
The group $OC(L)$ is called the coincidence isometry group (CIG) of $L$. For $Y\in OC(L)$, let
\be
\Sigma_{L}(Y) = [L: L\cap YL].
\ee 
If the structure matrix of $L$ is $A$, then the group $OC(L)$ is isomorphic to  $O_{n}(\mathbb{R}^{n})\cap (AGL_{n}(\mathbb{Q})A^{-1})$. If further $A$ is a rational matrix (in particular, this is the case if $L=\mathbb{Z}^{n}$), then 
\be
OC(L)= O_{n}(\mathbb{Q}):= O_{n}(\mathbb{R}^{n})\cap GL_{n}(\mathbb{Q}),
\ee
i.e. the corresponding group $OC(L)$ is formed by the rational orthogonal matrices. In this case, $OC(L)$ is generated by the reflections defined by the nonzero vectors of $L$ (see Zou (2006)). For $Y\in O_{n}(\mathbb{Q})$, write
\bea
Y = \frac{t}{q}Z,
\eea
where $t,q\in\mathbb{Z}_{+}$ such that $\gcd(t,q)=1$ and $Z$ is an integer matrix such that $\gcd(Z)=1$. Then since $\det Z \in \mathbb{Z}$ and 
\be
\det Y = (\frac{t}{q})^{n}\det Z = \pm 1,
\ee
we must have $t=1$ and
\bea
Y = \frac{1}{q}Z.
\eea
Let $q_{i}$ ($i=1,\ldots, n$) be the diagonal elements of the normal form of $Z$ (Fortes (1983)) and let
\bea
q_{(i)}=\frac{q}{\gcd(q,q_{i})}.
\eea
Then Fortes' result says that
\bea
\Sigma_{\mathbb{Z}^{n}}(Y)=q_{(1)}q_{(2)}\cdots q_{(n)}.
\eea
Since $\gcd(Z)=1$, $q_{1}=1$ and thus $q_{(1)}=q$. Furthermore, we have the following basic lemma.
\begin{lemma}
Let $Y$, $Z$, and $q$ be as in (2.3), and let $q_{i}$ ($i=1,\ldots, n$) be the diagonal elements of the normal form of $Z$. Then $q_{i}q_{n-i+1}=q^{2}$.
\end{lemma}
\begin{proof} From the discussion above, there are $P,Q\in GL_{n}(\mathbb{Z})$ (integer matrices with $\det=\pm 1$) such that
\bea
PYQ = \frac{1}{q}\left(\begin{array}{cccc} 
1 & {} & {} & {}\\
{} & q_{2} & {} & {}\\
{} & {} & \ddots & {}\\
{} & {} & {} & q_{n} \end{array}\right).
\eea
Taking inverses we have
\bea
Q^{-1}Y^{-1}P^{-1} &=& q\left(\begin{array}{cccc} 
1 & {} & {} & {}\\
{} & q_{2}^{-1} & {} & {}\\
{} & {} & \ddots & {}\\
{} & {} & {} & q_{n}^{-1} \end{array}\right)\nonumber \\
&=&
\frac{q}{q_{n}}\left(\begin{array}{cccc} 
q_{n} & {} & {} & {}\\
{} & q_{n}/q_{2} & {} & {}\\
{} & {} & \ddots & {}\\
{} & {} & {} & 1 \end{array}\right).
\eea
Note that the last integer matrix has the normal form
\bea
\left(\begin{array}{cccc} 
1 & {} & {} & {}\\
{} & q_{n}/q_{n-1} & {} & {}\\
{} & {} & \ddots & {}\\
{} & {} & {} & q_{n} \end{array}\right).
\eea
However, if we take transposes on both sides of (2.6), we have
\bea
Q^{T}Y^{T}P^{T} = \frac{1}{q}\left(\begin{array}{cccc} 
1 & {} & {} & {}\\
{} & q_{2} & {} & {}\\
{} & {} & \ddots & {}\\
{} & {} & {} & q_{n} \end{array}\right).
\eea
Since $Y^{-1}=Y^{T}$, by the uniqueness of the normal form, (2.7), (2.8), and (2.9) imply that $q_{n}=q^{2}$, which in turn implies the lemma.
\end{proof}
We will use this lemma to derive our index formulas in the next section.
\par
\section{Index formulas}
In this section, we assume $L=\mathbb{Z}^{n}$ and write $\Sigma(Y)$ for $\Sigma_{\mathbb{Z}^{n}}(Y)$. We begin with an immediate consequence of Lemma 2.1.
\begin{theorem}
Let $Y$, $Z$, $q$ be as in (2.3), and let $\delta_{i}$ be the greatest common divisor of the determinants of the $i\times i$ minors of $Z$ ($i=1,\ldots, n$). Then 
\bea
\Sigma(Y)=\frac{q^{m}}{\delta_{m}},
\eea
where $m=[n/2]$ is the integer part of $n/2$. 
\end{theorem}
\begin{proof}
Since the diagonal elements of the normal form of $Z$ satisfy $q_{i}\mid q_{i+1}$, by Lemma 2.1, $q_{i}\mid q$ if $i\le m$, and $q\mid q_{i}$ if $i>m$. Thus the $q_{(i)}$ defined in (2.4) are given by 
\be
q_{(i)}=\left\{ \begin{array}{rll} 
q/q_{i} & \mbox{if} & i\le m,\\
1 & \mbox{if} & i>m.
\end{array}\right.
\ee
Therefore, since the greatest common divisors of the determinants of the minors of $Z$ and its normal form are the same (see for example p. 458 and p. 485 in Artin (1991)), using (2.5) we have
\be
\Sigma(Y)=\frac{q^{m}}{q_{1}\cdots q_{m}}=\frac{q^{m}}{\delta_{m}}.
\ee 
\end{proof}
\par
{\bf Remark.} Note that since $\delta_{1}=1$, for $n\le 3$, formula (3.1) simplifies to the known result $\Sigma(Y)=q$ (see Baake (1997)). Note also that for $n=4$ and $5$, the formula is the same: $\Sigma(Y)=q^{2}/\delta_{2}$.
\par
Since the group $OC(\mathbb{Z}^{n})$ is generated by reflections defined by the nonzero vectors of $\mathbb{Z}^{n}$, we now turn to the reflections. Let
\bea
0\ne v=\sum_{i}^{n}a_{i}e_{i},\quad a_{i}\in\mathbb{Z},\quad 1\le i\le n.
\eea
Since we are interested in the reflection defined by $v$, we can always assume that $\gcd(a_{1},\ldots, a_{n})=1$. Under the canonical basis, the reflection of $\mathbb{R}^{n}$ defined by $v$ has the matrix
\bea
R_{v}=I-2\frac{vv^{T}}{v^{T}v}=\frac{1}{v^{T}v}(v^{T}vI-2vv^{T}),
\eea
where $v^{T}$ is the transpose of $v$.
\par
\begin{lemma}
Let $v$ be as above. Then
\bea
\gcd(v^{T}vI-2vv^{T})=\left\{ \begin{array}{rlll} 
1 & \mbox{if} & v^{T}v & \mbox{is odd},\\
2 & \mbox{if} & v^{T}v & \mbox{is even}.
\end{array}\right.
\eea
\end{lemma}
\begin{proof}
The $i$-th row of the matrix $v^{T}vI-2vv^{T}$ is 
\be
r_{i}=(-2a_{i}a_{1},\ldots,v^{T}v-2a_{i}^{2},\ldots,-2a_{i}a_{n}).
\ee 
Let $d_{i}=\gcd(r_{i})$ or $\gcd(r_{i}/2)$ according to whether $v^{T}v$ is odd or even. We claim that if a prime $p$ divides $d_{i}$, then it divides $a_{i}$. In fact, if 
\be
t_{i}=\gcd(a_{1},\ldots,a_{i-1},a_{i+1},\ldots,a_{n}), 
\ee
then
\be
d_{i}=\left\{ \begin{array}{rlll} 
\gcd(2a_{i}t_{i},v^{T}v-2a_{i}^{2}) & \mbox{if} & v^{T}v & \mbox{is odd},\\
\gcd(a_{i}t_{i},(v^{T}v-2a_{i}^{2})/2) & \mbox{if} & v^{T}v & \mbox{is even}.
\end{array}\right.
\ee
So if $p\mid d_{i}$ but $p\nmid a_{i}$, then $p\mid t_{i}$, implies that $p\mid\sum_{k\ne i}^{n}a_{k}^{2}$. However, this would imply that $p$ also divides 
\be
a_{i}^{2}=2a_{i}^{2}-v^{T}v+\sum_{k\ne i}^{n}a_{k}^{2}, 
\ee
which is a contradiction. Now the lemma follows since by our assumption $\gcd(a_{1},\ldots,a_{n})=1$.
\end{proof}
\par
It follows from this lemma that if we write 
\bea
R_{v}=\frac{1}{q}T_{v}
\eea
as in (2.3), then $q=v^{T}v$ or $v^{T}v/2$ depending on whether $v^{T}v$ is odd or even. Moreover, we have the following lemma: 
\begin{lemma}
Assume that $n>2$. If $q_{1},q_{2},\ldots,q_{n}$ are the diagonal elements of the normal form of $T_{v}$, then $q_{2}=q$.
\end{lemma}
{\bf Remark.} Note that for $n=2$ we have $q_{2}=q^{2}$ by Lemma 2.1. Note also that it follows from this lemma that $q_{2}=\cdots=q_{n-1}=q$. 
\par
\begin{proof}
To prove the lemma, recall that $q_{i}=\delta_{i}/\delta_{i-1}$, where $\delta_{i}$ is the greatest common divisor of the determinants of all $i\times i$ minors of $T_{v}$ (see for example p. 458 and p. 485 in Artin (1991)). Since $\delta_{1}=1$, we need to prove $\delta_{2}=q$. We give the detail for the case that $v^{T}v$ is odd, since it will be clear from the discussion that the same argument works for the even case. Consider the $2\times 2$ minors of $T_{v}$. If a $2\times 2$ minor $M$ does not involve any diagonal element, then $\det(M)=0$. If it involves diagonal element(s), then there are basically two possibilities:
\be
\left(\begin{array}{cc} 
v^{T}v-2a_{i}^{2} & {}-2a_{i}a_{j}\\
{}-2a_{t}a_{i} & {}-2a_{t}a_{j}
\end{array}\right) \quad \mbox{or} \quad 
\left(\begin{array}{cc} 
v^{T}v-2a_{i}^{2} & {}-2a_{i}a_{j}\\
{}-2a_{j}a_{i} & v^{T}v-2a_{j}^{2}
\end{array}\right). 
\ee
It is clear that the determinants of both matrices have the factor $v^{T}v$, hence our claim follows.
\end{proof} 
\par
We now give a formula for the index of a reflection.
\begin{theorem}
Let
\be
0\ne v=\sum_{i=1}^{n}a_{i}e_{i}\in\mathbb{Z}^{n}\quad\mbox{with}\quad \gcd(a_{1},\ldots,a_{n})=1.
\ee
Then
\bea
\Sigma(R_{v})=\left\{ \begin{array}{rlll} 
v^{T}v & \mbox{if} & v^{T}v & \mbox{is odd},\\
v^{T}v/2 & \mbox{if} & v^{T}v & \mbox{is even}.
\end{array}\right.
\eea
\end{theorem} 
\begin{proof}
This is an immediate consequence of Lemma 2.1, Lemma 3.1, and Lemma 3.2.
\end{proof}
\par
This theorem provides the base for using induction to obtain some interesting results. As an example, we will prove a proposition.
\par
For $i=1,\ldots,k$, let
\be
0\ne v_{i}=\sum_{j=1}^{n}a_{ji}e_{j}\in\mathbb{Z}^{n}\quad\mbox{with}\quad \gcd(a_{1i},\ldots,a_{ni})=1.
\ee
Define the integers $r_{i}$ to be $v_{i}^{T}v_{i}$ or $v_{i}^{T}v_{i}/2$ depending on whether $v_{i}^{T}v_{i}$ is odd or even. 
\begin{proposition} 
Assume that $\gcd(r_{i},r_{j})=1$ for $i\ne j$. Then
\be
\Sigma(R_{v_{1}}\cdots R_{v_{k}})=r_{1}\cdots r_{k}.
\ee  
\end{proposition}
\par
This proposition follows from Theorem 3.2 and the following lemma:
\begin{lemma}
Let $R_{i}\in O_{n}(\mathbb{Q})$ ($i=1,2$) be reflections and write $R_{i}=(1/r_{i})S_{i}$ as in (2.3). Assume that $\gcd(r_{1},r_{2})=1$ and the normal forms of $S_{i}$ are
\be
\left(\begin{array}{cccc} 
1 & {} & {} & {}\\
{} & r_{i} & {} & {}\\
{} & {} & \ddots & {}\\
{} & {} & {} & r_{i}^{2} \end{array}\right), \quad i=1,2.
\ee
If we write $R_{1}R_{2}=(1/r)R$ as in (2.3), then $r=r_{1}r_{2}$ and the matrix $R$ has the normal form
\be
\left(\begin{array}{cccc} 
1 & {} & {} & {}\\
{} & r & {} & {}\\
{} & {} & \ddots & {}\\
{} & {} & {} & r^{2} \end{array}\right).
\ee
\end{lemma}
\begin{proof}
Let the normal form of $R$ be
\be
\left(\begin{array}{cccc} 
d_{1} & {} & {} & {}\\
{} & d_{2} & {} & {}\\
{} & {} & \ddots & {}\\
{} & {} & {} & d_{n} \end{array}\right).
\ee
Then we know that $d_{1}=1$ and $d_{n}=r^{2}$, so it remains to prove that $d_{2}=r$. Consider
\be
R_{1}R_{2}=\frac{1}{r_{1}r_{2}}S_{1}S_{2}.
\ee
If the normal form of $S_{1}S_{2}$ is
\be
\left(\begin{array}{cccc} 
c_{1} & {} & {} & {}\\
{} & c_{2} & {} & {}\\
{} & {} & \ddots & {}\\
{} & {} & {} & c_{n} \end{array}\right),
\ee
then $r=r_{1}r_{2}/c_{1}$ and $d_{i}=c_{i}/c_{1}$. For an integer matrix $A$, denote by $\delta_{i}(A)$ the greatest common divisor of the determinants of the $i\times i$ minors of $A$. Then $\delta_{i}(S_{1}S_{2})$ are identical to those of (compare with the proof of Lemma 2.1) 
\bea
R^{\prime}=\left(\begin{array}{cccc} 
1 & {} & {} & {}\\
{} & r_{1} & {} & {}\\
{} & {} & \ddots & {}\\
{} & {} & {} & r_{1}^{2} \end{array}\right)S_{2}^{\prime},
\eea
where $S_{2}^{\prime}$ is obtained from $S_{2}$ by left multiplying by an element from $GL_{n}(\mathbb{Z})$. Therefore $\gcd(S_{2}^{\prime})=1$ and $\delta_{i}(S_{2}^{\prime})=\delta_{i}(S_{2})$. Thus if a prime $p$ divides $c_{1}=\gcd(R^{\prime})$, it must divide $r_{1}$. Similarly, $p$ also divides $r_{2}$. But $\gcd(r_{1},r_{2})=1$, so $\gcd(R^{\prime})=c_{1}=1$. Thus $r=r_{1}r_{2}$ and $d_{2}=c_{2}$.
\par
Now $c_{2}=\delta_{2}(R^{\prime})/c_{1}=\delta_{2}(R^{\prime})$. Since every $2\times 2$ minor of $R^{\prime}$ contains at least one row of $S_{2}^{\prime}$ multiplied by $r_{1}$ or $r_{1}^{2}$, we see that $r_{1}r_{2}=r_{1}\delta_{2}(S_{2}^{\prime})\mid \delta_{2}(R^{\prime})$, implies that $r\mid d_{2}$. But $d_{2}\le r$ by Lemma 2.1, so $d_{2}=r$.
\end{proof}
\par
It should be pointed out that without the assumption that $\gcd(r_{i},r_{j})=1$, the result of Proposition 3.1 does not hold. This can be seen by noting that the square of a reflection is the identity.
\section{Concluding remarks}
\par
It is known that certain positive integers can not be the coincidence isometry indices for the lattice $\mathbb{Z}^{n}$ when $n\le 4$. In particular, it is well-known that in dimension 3, the indices assume precisely the odd positive integers (Grimmer (1973)). One may ask what happens when $n\ge 4$. This question can be answered by using the results of the present work together with some known facts about the square sums of integers. Recall that a theorem due to Legendre and Gauss says that a positive integer can be expressed as a sum of three squares if and only if it is not the form $4^{m}(8k+7)$ (Adler and Coury (1995), p. 236). It follows from this theorem that every {\it odd} positive integer can be written as a sum of four integers with $\gcd = 1$. To see this, note that for $n\ge 0$, if
\be
4n+1=a^{2}+b^{2}+c^{2},
\ee
then since every square is congruent to $0$ or $1$ modulo 4, exactly one of $a$, $b$, and $c$ is odd. Assume that $c$ is odd and write
\be
a=2u,\quad b=2v, \quad c=2t+1,
\ee
then 
\be
(u+v)^{2}+(u-v)^{2}+t^{2}+(t+1)^{2}=2n+1.
\ee
Thus, by Theorem 3.2, the indices provided by $OC(\mathbb{Z}^{n})$ ($n\le 4$) cover all the positive odd integers, though these indices miss the powers of $2^{k}$ for $k>1$ (see Baake (1997)). However, for $n=5$, all the positive integers are covered. To see this, note that since $2^{k}-1$ is odd, from the above discussion, $2^{k}$ can be expressed as a sum of five squares with $\gcd=1$, so Theorem 3.2 and Proposition 3.1 imply the result.
\par
The index formulas for the coincidence isometries of the lattice $\mathbb{Z}^{n}$ provided in this work are quite explicit. However, the computations are more involved in the general cases, and one should not expect to have formulas as explicit. The connection between the coincidence index formulas and the related formulas in number theory deserves further attention (see Baake (1997)).     
\medskip

\end{document}